\documentclass[12pt]{article}
\usepackage{amsmath,amssymb, amsfonts, amsthm,amscd}
\usepackage[T2A]{fontenc}
\usepackage[cp1251]{inputenc}
\usepackage[russian,ukrainian,english]{babel}
\usepackage{geometry}

\newcommand{\case}[1]{\textbf{Case $\mathbf{#1}$.}\ }

\newtheorem{Theorem}{Theorem}
\newtheorem{Corollary}{Corollary}
\newtheorem{Remark}{Remark}
\newtheorem{Lemma}{Lemma}
\newtheorem{Proposition}{Proposition}

\newenvironment{LemmaProof}{\textbf{Proof. }}{\par\noindent\textbf{The Lemma is proved.}}

\title{{\Large \textbf{On the extremal values of the number of vertices with an interval spectrum on the set of
                proper edge colorings of the graph of the $n$-dimensional cube}}}

\author{\normalsize A.M. Khachatryan$^1$, R.R. Kamalian$^2$}

\date{}

\begin{document}

\maketitle

$\\^1$Ijevan Branch of Yerevan State University, e-mail:
khachatryanarpine@gmail.com $\\^2$The Institute for Informatics and
Automation Problems of NAS RA, \\e-mail: rrkamalian@yahoo.com
\bigskip \bigskip

\begin{abstract}
For an undirected, simple, finite, connected graph $G$, we denote by
$V(G)$ and $E(G)$ the sets of its vertices and edges, respectively.
A function $\varphi:E(G)\rightarrow \{1,...,t\}$ is called a proper
edge $t$-coloring of a graph $G$, if adjacent edges are colored
differently and each of $t$ colors is used. The least value of $t$
for which there exists a proper edge $t$-coloring of a graph $G$ is
denoted by $\chi'(G)$. For any graph $G$, and for any integer $t$
satisfying the inequality $\chi'(G)\leq t\leq |E(G)|$, we denote by
$\alpha(G,t)$ the set of all proper edge $t$-colorings of $G$. Let
us also define a set $\alpha(G)$ of all proper edge colorings of a
graph $G$:
$$
\alpha(G)\equiv\bigcup_{t=\chi'(G)}^{|E(G)|}\alpha(G,t).
$$

An arbitrary nonempty finite subset of consecutive integers is
called an interval. If $\varphi\in\alpha(G)$ and $x\in V(G)$, then
the set of colors of edges of $G$ which are incident with $x$ is
denoted by $S_G(x,\varphi)$ and is called a spectrum of the vertex
$x$ of the graph $G$ at the proper edge coloring $\varphi$. If $G$
is a graph and $\varphi\in\alpha(G)$, then define
$f_G(\varphi)\equiv|\{x\in V(G)/S_G(x,\varphi) \textrm{ is an
interval}\}|$.

For a graph $G$ and any integer $t$, satisfying the inequality
$\chi'(G)\leq t\leq |E(G)|$, we define:
$$
\mu_1(G,t)\equiv\min_{\varphi\in\alpha(G,t)}f_G(\varphi),\qquad
\mu_2(G,t)\equiv\max_{\varphi\in\alpha(G,t)}f_G(\varphi).
$$

For any graph $G$, we set:
$$
\mu_{11}(G)\equiv\min_{\chi'(G)\leq t\leq|E(G)|}\mu_1(G,t),\qquad
\mu_{12}(G)\equiv\max_{\chi'(G)\leq t\leq|E(G)|}\mu_1(G,t),
$$
$$
\mu_{21}(G)\equiv\min_{\chi'(G)\leq t\leq|E(G)|}\mu_2(G,t),\qquad
\mu_{22}(G)\equiv\max_{\chi'(G)\leq t\leq|E(G)|}\mu_2(G,t).
$$

For any positive integer $n$, the exact values of the parameters
$\mu_{11}$, $\mu_{12}$, $\mu_{21}$ and $\mu_{22}$ are found for the
graph of the $n$-dimensional cube.


Math. Classification: 05C15
\end{abstract}

We consider finite, undirected, connected graphs without loops and
multiple edges containing at least one edge. For any graph $G$, we
denote by $V(G)$ and $E(G)$ the sets of vertices and edges of $G$,
respectively. For any $x\in V(G)$, $d_G(x)$ denotes the degree of
the vertex $x$ in $G$. For a graph $G$, $\delta(G)$ and $\Delta(G)$
denote the minimum and maximum degrees of vertices in $G$,
respectively. For a graph $G$, and for any $V_0\subseteq V(G)$, we
denote by $G[V_0]$ the subgraph of the graph $G$ induced
\cite{West1} by the subset $V_0$ of its vertices. For arbitrary
graphs $G$ and $H$, $G\Box H$ denotes their cartesian product
\cite{West1}.

For any positive integer $n$, we denote by $Q_n$ the graph of the
$n$-dimensional cube \cite{Furedi}. Clearly, for any positive
integer $n$, $Q_n$ is a regular bipartite graph with $|V(Q_n)|=2^n$,
$|E(Q_n)|=n\cdot 2^{n-1}$, $\Delta(Q_n)=n$.

An arbitrary nonempty finite subset of consecutive integers is
called an interval. An interval with the minimum element $p$ and the
maximum element $q$ is denoted by $[p,q]$.

A function $\varphi:E(G)\rightarrow [1,t]$ is called a proper edge
$t$-coloring of a graph $G$, if each of $t$ colors is used, and
adjacent edges are colored differently.

The minimum value of $t$ for which there exists a proper edge
$t$-coloring of a graph $G$ is denoted by $\chi'(G)$ \cite{Vizing2}.

For any graph $G$, and for any $t\in[\chi'(G),|E(G)|]$, we denote by
$\alpha(G,t)$ the set of all proper edge $t$-colorings of $G$.

Let us also define a set $\alpha(G)$ of all proper edge colorings of
a graph $G$:
$$
\alpha(G)\equiv\bigcup_{t=\chi'(G)}^{|E(G)|}\alpha(G,t).
$$

If $\varphi\in\alpha(G)$ and $x\in V(G)$, then the set
$\{\varphi(e)/ e\in E(G), e \textrm{ is incident with } x$\} is
called a spectrum of the vertex $x$ of the graph $G$ at the proper
edge coloring $\varphi$ and is denoted by $S_G(x,\varphi)$; if
$S_G(x,\varphi)$ is an interval, we say that $\varphi$ is interval
in $x$.

If $G$ is a graph, $\varphi\in\alpha(G)$, $R\subseteq V(G)$, then we
say, that $\varphi$ is interval on $R$ iff for $\forall x\in R$,
$\varphi$ is interval in $x$. We say, that a subset $R$ of vertices
of a graph $G$ has an $i$-property iff there exists
$\varphi\in\alpha(G)$ interval on $R$. If $G$ is a graph, and a
subset $R$ of its vertices has an $i$-property, we denote by
$w_R(G)$ and $W_R(G)$ (omiting the index in these notations in the
peculiar case of $R=V(G)$) the minimum and the maximum value of $t$,
respectively, for which $\exists\varphi\in\alpha(G,t)$ interval on
$R$ \cite{Oranj3, Asratian4, Diss5} (see also \cite{Vestnik15}).

If $G$ is a graph, $\varphi\in\alpha(G)$, then set
$V_{int}(G,\varphi)\equiv\{x\in V(G)/S_G(x,\varphi) \textrm{ is an
interval}\}$ and $f_G(\varphi)\equiv|V_{int}(G,\varphi)|$. A proper
edge coloring $\varphi\in\alpha(G)$ is called an interval edge
coloring \cite{Oranj3, Asratian4, Diss5} of the graph $G$ iff
$f_G(\varphi)=|V(G)|$. The set of all graphs having an interval edge
coloring is denoted by $\mathfrak{N}$. The simplest example of the
graph which doesn't belong to $\mathfrak{N}$ is $K_3$. The terms and
concepts which are not defined can be found in \cite{West1}.

For a graph $G$, and for any $t\in[\chi'(G),|E(G)|]$, we set
\cite{Mebius6}:
$$
\mu_1(G,t)\equiv\min_{\varphi\in\alpha(G,t)}f_G(\varphi),\qquad
\mu_2(G,t)\equiv\max_{\varphi\in\alpha(G,t)}f_G(\varphi).
$$

For any graph $G$, we set \cite{Mebius6}:
$$
\mu_{11}(G)\equiv\min_{\chi'(G)\leq t\leq|E(G)|}\mu_1(G,t),\qquad
\mu_{12}(G)\equiv\max_{\chi'(G)\leq t\leq|E(G)|}\mu_1(G,t),
$$
$$
\mu_{21}(G)\equiv\min_{\chi'(G)\leq t\leq|E(G)|}\mu_2(G,t),\qquad
\mu_{22}(G)\equiv\max_{\chi'(G)\leq t\leq|E(G)|}\mu_2(G,t).
$$

Clearly, the parameters $\mu_{11}$, $\mu_{12}$, $\mu_{21}$ and
$\mu_{22}$ are correctly defined for an arbitrary graph.

Let us note that exact values of the parameters $\mu_{12}$ and
$\mu_{21}$ have certain game interpretations. Suppose that all edges
of a graph $G$ are colored in the game of Alice and Bob with
antagonistic interests and asymmetric distribution of roles. Alice
determines the number $t$ of colors in the future proper edge
coloring $\varphi$ of the graph $G$, satisfying the condition
$t\in[\chi'(G),|E(G)|]$, Bob colors edges of $G$ with $t$ colors.

When Alice aspires to maximize, Bob aspires to minimize the value of
the function $f_G(\varphi)$, and both players choose their best
strategies, then at the finish of the game exactly $\mu_{12}(G)$
vertices of $G$ will receive an interval spectrum.

When Alice aspires to minimize, Bob aspires to maximize the value of
the function $f_G(\varphi)$, and both players choose their best
strategies, then at the finish of the game exactly $\mu_{21}(G)$
vertices of $G$ will receive an interval spectrum.

The exact values of the parameters $\mu_{11}$, $\mu_{12}$,
$\mu_{21}$ and $\mu_{22}$ are found for simple paths, simple cycles
and simple cycles with a chord \cite{Simple7, Akunq}, "M\"{o}bius
ladders" \cite{Mebius6, Minchev}, complete graphs \cite{Arpine8},
complete bipartite graphs \cite{Arpine9, Arpine10}, prisms
\cite{Arpine11, Minchev}, $n$-dimensional cubes \cite{Arpine11,
Nikolaev12} and the Petersen graph \cite{Petersen}. The exact values
of $\mu_{11}$ and $\mu_{22}$ for trees are found in \cite{Evg13}.
The exact value of $\mu_{12}$ for an arbitrary tree is found in
\cite{Trees14} (see also \cite{Algorithm, Tree_Kontr}).

In this paper, for any positive integer $n$, we determine the exact
values of the parameters $\mu_{11}$, $\mu_{12}$, $\mu_{21}$ and
$\mu_{22}$ for the graph $Q_n$.

First we recall some known results.

\begin{Lemma}\cite{Oranj3, Asratian4, Diss5}\label{Thm1}
If $G\in\mathfrak{N}$, then $\chi'(G)=\Delta(G)$.
\end{Lemma}

\begin{Corollary}\label{cor1}
If $G\in\mathfrak{N}$, then $\Delta(G)=\chi'(G)\leq w(G)\leq
W(G)\leq|E(G)|$.
\end{Corollary}

\begin{Corollary}\label{cor2}
If $G$ is a regular graph, then $G\in\mathfrak{N}$ iff
$\chi'(G)=\Delta(G)$.
\end{Corollary}

\begin{Corollary}\label{cor3}
If $G$ is a regular graph with $\chi'(G)=\Delta(G)$, then
$G\in\mathfrak{N}$ and $w(G)=\chi'(G)=\Delta(G)$.
\end{Corollary}

\begin{Corollary}\label{cor3n}
If $G$ is a regular bipartite graph, then $G\in\mathfrak{N}$ and
$w(G)=\chi'(G)=\Delta(G)$.
\end{Corollary}

\begin{Corollary}\label{cor4}
If $G\in\mathfrak{N}$, and a subset $R\subseteq V(G)$ has an
$i$-property, then $\Delta(G)=\chi'(G)\leq  w_R(G)\leq w(G)\leq
W(G)\leq W_R(G)\leq |E(G)|$.
\end{Corollary}

\begin{Lemma}\cite{Oranj3, Asratian4, Diss5}\label{Thm2}
If $R$ is the set of all vertices of an arbitrary part of a
bipartite graph $G$, then:
\begin{enumerate}
  \item $R$ has an $i$-property,
  \item $W_R(G)=|E(G)|$,
  \item for any $t\in[w_R(G),W_R(G)]$, there exists $\varphi_t\in\alpha(G,t)$ interval on $R$.
\end{enumerate}
\end{Lemma}

\begin{Corollary}\label{cor5}
If $R$ is the set of all vertices of an arbitrary part of a regular
bipartite graph  $G$, then:
\begin{enumerate}
  \item $R$ has an $i$-property,
  \item $W_R(G)=|E(G)|$,
  \item $w_R(G)=\Delta(G)$,
  \item for any $t\in[w_R(G),W_R(G)]$, there exists $\varphi_t\in\alpha(G,t)$ interval on $R$.
\end{enumerate}
\end{Corollary}

\begin{Lemma}\cite{Luhansk17}\label{Thm3}
If $G$ is a graph with $\delta(G)\geq 2$,
$\varphi\in\alpha(G,|E(G)|)$, $V_{int}(G,\varphi)\neq\emptyset$,
then $G[V_{int}(G,\varphi)]$ is a forest each connected component of
which is a simple path.
\end{Lemma}

\begin{Lemma}\cite{Simple7, Akunq}\label{Thm4}
If $G$ is a regular graph with $\chi'(G)=\Delta(G)$, then
$\mu_{12}(G)=|V(G)|$.
\end{Lemma}

\begin{LemmaProof}
Let us note that for an arbitrary $\varphi\in\alpha(G,\Delta(G))$,
$f_G(\varphi)=|V(G)|$. Consequently,
$|V(G)|=\mu_1(G,\Delta(G))\leq\max_{t\in
[\chi'(G),|E(G)|]}\mu_1(G,t)=\mu_{12}(G)\leq|V(G)|$. Hence,
$\mu_{12}(G)=|V(G)|$.
\end{LemmaProof}

\begin{Lemma}\cite{Simple7, Akunq}\label{Thm5}
For an arbitrary graph $G$, the following inequalities hold:
$\mu_{11}(G)\leq\mu_{12}(G)\leq\mu_{22}(G)$, $\mu_{11}(G)\leq
\mu_{21}(G)\leq \mu_{22}(G)$.
\end{Lemma}

\begin{LemmaProof}
Clearly,
$$
\mu_{11}(G)=\min_{\chi'(G)\leq
t\leq|E(G)|}\mu_1(G,t)\leq\max_{\chi'(G)\leq
t\leq|E(G)|}\mu_1(G,t)=\mu_{12}(G).
$$
and
$$
\mu_{21}(G)=\min_{\chi'(G)\leq
t\leq|E(G)|}\mu_2(G,t)\leq\max_{\chi'(G)\leq
t\leq|E(G)|}\mu_2(G,t)=\mu_{22}(G).
$$

The inequality $\mu_{12}(G)\leq\mu_{22}(G)$ is provided by the
existence of an integer $t'\in[\chi'(G),|E(G)|]$, for which the
relation
$$
\mu_{12}(G)=\max_{\chi'(G)\leq
t\leq|E(G)|}\mu_1(G,t)=\mu_1(G,t')=\min_{\varphi\in\alpha(G,t')}f_G(\varphi)\leq
$$
$$
\leq\max_{\varphi\in\alpha(G,t')}f_G(\varphi)=\mu_2(G,t')\leq\max_{\chi'(G)\leq
t\leq|E(G)|}\mu_2(G,t)=\mu_{22}(G)
$$
is true.

The inequality $\mu_{11}(G)\leq\mu_{21}(G)$ is provided by the
existence of an integer $t''\in[\chi'(G),|E(G)|]$, for which the
relation
$$
\mu_{21}(G)=\min_{\chi'(G)\leq
t\leq|E(G)|}\mu_2(G,t)=\mu_2(G,t'')=\max_{\varphi\in\alpha(G,t'')}f_G(\varphi)\geq
$$
$$
\geq\min_{\varphi\in\alpha(G,t'')}f_G(\varphi)=\mu_1(G,t'')\geq\min_{\chi'(G)\leq
t\leq|E(G)|}\mu_1(G,t)=\mu_{11}(G)
$$
is true.

\end{LemmaProof}

\begin{Remark}\cite{Simple7, Akunq}
There are graphs $G$ with $\mu_{21}(G)<\mu_{12}(G)$. There are
graphs $G$ with $\mu_{12}(G)<\mu_{21}(G)$.
\end{Remark}

\begin{Lemma}\label{lem1}\cite{Zas, Archive}
If a subset $V_0$ of the set of vertices of the graph $Q_3$ contains
at least $6$ vertices, then at least one of the following two
statements is true:
\begin{enumerate}
  \item there exist such vertices $x_1, x_2, x_3, x_4$ in $V_0$ that $Q_3[\{x_1, x_2, x_3, x_4\}]\cong K_{3,1}$,
  \item there exist such vertices $y_1, y_2, y_3, y_4, y_5, y_6$ in $V_0$ that $Q_3[\{y_1, y_2, y_3, y_4, y_5, y_6\}]\cong C_6$.
\end{enumerate}
\end{Lemma}

\begin{Lemma}\label{Thm6}\cite{Zas, Archive}
If a subset $V_0$ of the set of vertices of the graph $Q_n$
$(n\geq4)$ contains at least $2^{n-1}+1$ vertices, then at least one
of the following two statements is true:
\begin{enumerate}
  \item there exist such vertices $x_1, x_2, x_3, x_4$ in $V_0$ that $Q_n[\{x_1, x_2, x_3, x_4\}]\cong K_{3,1}$,
  \item there exist such vertices $y_1, y_2, y_3, y_4, y_5, y_6, y_7, y_8$ in $V_0$ that $Q_n[\{y_1, y_2, y_3, y_4, y_5, y_6, y_7, y_8\}]\cong C_8$.
\end{enumerate}
\end{Lemma}

If $G$ is a graph with $\chi'(G)=\Delta(G)$,
$t\in[\Delta(G),|E(G)|]$, $\xi\in\alpha(G,t)$, then for any
$j\in[1,t]$, we denote by $E(G,\xi,j)$ the set of all edges of $G$
colored by the color $j$ at the proper edge $t$-coloring $\xi$;
$\xi$ is called a harmonic \cite{Arpine9} edge $t$-coloring of the
graph $G$, if for any $i\in[1,\Delta(G)]$, the set
$$
\bigcup_{1\leq j\leq t,\;j\equiv i(mod(\Delta(G)))}E(G,\xi,j)
$$
is a matching in $G$.

Suppose that $G$ is a graph with $\chi'(G)=\Delta(G)$, $t\in
[\Delta(G),|E(G)|]$, $\xi$ is a harmonic edge $t$-coloring of $G$.
Let us define a sequence $\xi^*_0,\xi^*_1,\ldots,\xi^*_{t-\chi'(G)}$
of proper edge colorings of the graph $G$.

Set $\xi^*_0\equiv\xi$.

\case{1} $t=\chi'(G)$. The desired sequence mentioned above is
already constructed.

\case{2} $t\in[\chi'(G)+1,|E(G)|]$. Suppose that
$j\in[1,t-\chi'(G)]$, and the proper edge colorings
$\xi^*_0,\ldots,\xi^*_{j-1}$ of the graph $G$ are already
constructed. Let us define $\xi^*_j$. For an arbitrary $e\in E(G)$,
set:
$$
\xi^*_j(e)\equiv\left\{
\begin{array}{ll}
\xi^*_{j-1}(e), & \textrm{if $\xi^*_{j-1}(e)\neq \max(\{\xi^*_{j-1}(e)/e\in E(G)\})$}\\
\xi^*_{j-1}(e)-\Delta(G), & \textrm{if $\xi^*_{j-1}(e)=
\max(\{\xi^*_{j-1}(e)/e\in E(G)\})$}.
\end{array}
\right.
$$

\begin{Remark}\label{Rem1}
Suppose that $G$ is a graph with $\chi'(G)=\Delta(G)$,
$t\in[\Delta(G),|E(G)|]$, $\xi$ is a harmonic edge $t$-coloring of
$G$. All proper edge colorings
$\xi^*_0,\xi^*_1,\ldots,\xi^*_{t-\chi'(G)}$ of the graph $G$ are
determined definitely.
\end{Remark}

\begin{Remark}\label{Rem2}
Suppose that $G$ is a graph with $\chi'(G)=\Delta(G)<|E(G)|$,
$t\in[1+\chi'(G),|E(G)|]$, $\xi$ is a harmonic edge $t$-coloring of
$G$. It is not difficult to see, that for any $j\in[1,t-\chi'(G)]$,
$\xi^*_j$ is a harmonic edge $(t-j)$-coloring of the graph $G$.
\end{Remark}

\begin{Remark}\label{Rem3}
Suppose that $G$ is a graph with $\chi'(G)=\Delta(G)<|E(G)|$,
$t\in[1+\chi'(G),|E(G)|]$, $\xi$ is a harmonic edge $t$-coloring of
$G$. Assume that $\xi$ is interval in some vertex $z_0\in V(G)$ with
$d_G(z_0)=\Delta(G)$. Then, for any $j\in[1,t-\chi'(G)]$, $\xi^*_j$
is interval in $z_0$.
\end{Remark}

From results of \cite{Simple7, Akunq} we have

\begin{Proposition}\label{lem2}
$\mu_{11}(Q_1)=2$, $\mu_{11}(Q_2)=1$, $\mu_{21}(Q_1)=2$,
$\mu_{21}(Q_2)=3$.
\end{Proposition}

\begin{Lemma}\label{lem3}
If an $r$-regular $(r\geq3)$ graph $G$ has $\varphi\in\alpha(G,r+1)$
with $f_G(\varphi)=0$, then the graph $G\Box K_2$ has
$\psi\in\alpha(G\Box K_2, r+2)$ with $f_{G\Box K_2}(\psi)=0$.
\end{Lemma}
\textbf{Proof} is evident.

\begin{Lemma}\label{lem4}
There exists $\varphi\in\alpha(Q_3,4)$ with $f_{Q_3}(\varphi)=0$.
There exists $\psi\in\alpha(Q_3,12)$ with $f_{Q_3}(\psi)=5$.
\end{Lemma}

\begin{LemmaProof}
Let $V(Q_3)=\{x_1,x_2,x_3,x_4,y_1,y_2,y_3,y_4\}$,
$E(Q_3)=\{(x_1,x_2), (x_2,x_3), (x_3,x_4), (x_1,x_4), \\(y_1,y_2),
(y_2,y_3), (y_3,y_4), (y_1,y_4), (x_1,y_1), (x_2,y_2), (x_3,y_3),
(x_4,y_4)\}$.

Set:
$\varphi((x_1,x_2))=\varphi((x_3,x_4))=\varphi((y_1,y_2))=\varphi((y_3,y_4))=1$,
$\varphi((x_1,x_4))=\varphi((x_2,x_3))=2$,
$\varphi((y_1,y_4))=\varphi((y_2,y_3))=3$,
$\varphi((x_1,y_1))=\varphi((x_2,y_2))=\varphi((x_3,y_3))=\varphi((x_4,y_4))=4$.

It is not difficult to see that $\varphi\in\alpha(Q_3,4)$ and
$f_{Q_3}(\varphi)=0$.

Set: $\psi((x_1,y_1))=1$, $\psi((y_1,y_2))=2$, $\psi((y_1,y_4))=3$,
$\psi((y_3,y_4))=4$, $\psi((x_4,y_4))=5$, $\psi((x_1,x_4))=6$,
$\psi((x_3,x_4))=7$, $\psi((x_3,y_3))=8$, $\psi((x_2,x_3))=9$,
$\psi((x_2,y_2))=10$, $\psi((x_1,x_2))=11$, $\psi((y_2,y_3))=12$.

It is not difficult to see that $\psi\in\alpha(Q_3,12)$ and
$f_{Q_3}(\psi)=5$.
\end{LemmaProof}

From Lemmas \ref{lem3} and \ref{lem4} we obtain

\begin{Proposition}\label{cor6}
For any integer $n\geq 3$, $\mu_{11}(Q_n)=0$.
\end{Proposition}

From Propositions \ref{lem2} and \ref{cor6} we obtain

\begin{Proposition}\label{Prop1}
For any positive integer $n$, $\mu_{11}(Q_n)=3-\min\{3,n\}$.
\end{Proposition}

From Lemma \ref{Thm4} we obtain

\begin{Proposition}\label{Prop2}
For any positive integer $n$, $\mu_{12}(Q_n)=2^n$.
\end{Proposition}

From Proposition \ref{Prop2} and Lemma \ref{Thm5} we obtain

\begin{Proposition}\label{Prop3}
For any positive integer $n$, $\mu_{22}(Q_n)=2^n$.
\end{Proposition}

\begin{Lemma}\label{lem5}
For any integer $n\geq4$, $\mu_{2}(Q_n,n\cdot 2^{n-1})\geq2^{n-1}$.
\end{Lemma}

\textbf{Proof} follows from Lemma \ref{Thm2}.

\begin{Lemma}\label{lem6}
For any integer $n\geq4$, $\mu_{2}(Q_n,n\cdot 2^{n-1})\leq2^{n-1}$.
\end{Lemma}

\begin{LemmaProof}
Assume the contrary. Then there exist an integer $n_0\geq4$ and
$\varphi_0\in\alpha(Q_{n_0},n_0\cdot 2^{n_0-1})$, for which
$f_{Q_{n_0}}(\varphi_0)\geq2^{n_0-1}+1$. By Lemma \ref{Thm3},
$Q_{n_0}[V_{int}(Q_{n_0},\varphi_0)]$ is a forest, each connected
component of which is a simple path. But it is incompatible with
Lemma \ref{Thm6}.
\end{LemmaProof}

From Lemmas \ref{lem5} and \ref{lem6} we obtain

\begin{Lemma}\label{lem7}
For any integer $n\geq4$, $\mu_{2}(Q_n,n\cdot 2^{n-1})=2^{n-1}$.
\end{Lemma}

\begin{Lemma}\label{lem8}
For arbitrary integers $n$ and $t$, satisfying the conditions
$n\geq4$, $n\leq t\leq n\cdot 2^{n-1}$, the inequality
$\mu_{2}(Q_n,t)\geq 2^{n-1}$ is true.
\end{Lemma}

\begin{LemmaProof}
Choose an arbitrary integer $n_0\geq 4$. It follows from Corollary
\ref{cor5}, that for any integer $t$, satisfying the inequality
$n_0\leq t\leq n_0\cdot 2^{n_0-1}$, there exists $\varphi_t\in
\alpha(Q_{n_0},t)$ with $f_{Q_{n_0}}(\varphi_t)\geq2^{n_0-1}$. It
means that for any integer $t$, satisfying the inequality $n_0\leq
t\leq n_0\cdot 2^{n_0-1}$, we also have $\mu_{2}(Q_{n_0},t)\geq
2^{n_0-1}$.
\end{LemmaProof}

From Lemmas \ref{lem7} and \ref{lem8} we obtain

\begin{Proposition}\label{cor7}
For any integer $n\geq4$, $\mu_{21}(Q_n)=2^{n-1}$.
\end{Proposition}

From Lemma \ref{lem4} we obtain

\begin{Lemma}\label{cor8}
$\mu_{2}(Q_3,12)\geq 5$.
\end{Lemma}

\begin{Lemma}\label{lem9}
$\mu_{2}(Q_3,12)\leq 5$.
\end{Lemma}

\begin{LemmaProof}
Assume the contrary: $\mu_{2}(Q_3,12)\geq 6$. It means that there
exists $\xi\in\alpha(Q_3,12)$ with $f_{Q_3}(\xi)\geq6$. By Lemma
\ref{Thm3}, $Q_3[V_{int}(Q_3,\xi)]$ is a forest each connected
component of which is a simple path. But it is incompatible with
Lemma \ref{lem1}.
\end{LemmaProof}

From Lemmas \ref{cor8} and \ref{lem9} we obtain

\begin{Lemma}\label{lem10}
$\mu_{2}(Q_3,12)=5$.
\end{Lemma}

\begin{Lemma}\label{lem11}
For any integer $t$, satisfying the condition $3\leq t\leq 12$, the
inequality $\mu_2(Q_3,t)\geq 5$ is true.
\end{Lemma}

\begin{LemmaProof}
It is not difficult to check that the proper edge $12$-coloring
$\psi$ of the graph $Q_3$ constructed for the proof of Lemma
\ref{lem4} is a harmonic edge $12$-coloring of $Q_3$. From Remark
\ref{Rem3} it follows that for any integer $t$ satisfying the
inequality $3\leq t\leq 11$, there exists $\psi_t\in\alpha(Q_3,t)$
with $V_{int}(Q_3,\psi)\subseteq V_{int}(Q_3,\psi_t)$. Hence, for
the same values of $t$, the inequality $\mu_2(Q_3,t)\geq
f_{Q_3}(\psi)=5$ is true.
\end{LemmaProof}

From Lemmas \ref{lem10} and \ref{lem11} we obtain

\begin{Proposition}\label{cor9}
$\mu_{21}(Q_3)=5$.
\end{Proposition}

From Propositions \ref{lem2}, \ref{cor7} and \ref{cor9} we obtain

\begin{Proposition}\label{Prop4}
For any positive integer $n$,
$$
\mu_{21}(Q_n)=2^{n-1}+\bigg\lceil1-\frac{\min\{4,n\}}{4}\bigg\rceil.
$$
\end{Proposition}

From Propositions \ref{Prop1} - \ref{Prop3} and \ref{Prop4} we
obtain

\begin{Theorem}\label{Thm7}\cite{Arpine11, Nikolaev12}
For any positive integer $n$, the following equalities hold:
$$
\mu_{11}(Q_n)=3-\min\{3,n\}, \qquad \mu_{12}(Q_n)=\mu_{22}(Q_n)=2^n,
$$
$$
\mu_{21}(Q_n)=2^{n-1}+\bigg\lceil1-\frac{\min\{4,n\}}{4}\bigg\rceil.
$$
\end{Theorem}

\end{document}